# From Morse-Smale to all Knots and Links

Robert Ghrist
Department of Mathematics
University of Texas, Austin
Austin, TX 78712
ghrist@math.utexas.edu

Todd Young
Department of Mathematics
Northwestern University
Evanston, IL 60208
young@math.ohiou.edu

**Abstract**

We analyse the topological (knot-theoretic) features of a certain codimension-one bifurcation of a partially hyperbolic fixed point in a flow on $\mathbb{R}^3$ originally described by Shil'nikov. By modifying how the invariant manifolds wrap around themselves, or "pleat," we may apply the theory of templates, or branched two-manifolds, to capture the topology of the flow. This analysis yields a class of flows which bifurcate from a Morse-Smale flow to a Smale flow containing periodic orbits of all knot and link types.

The goal of this paper is to show that flows in $\mathbb{R}^3$ with the most complex knot structure possible are accessible from the set of Morse-Smale flows via a comdimension one bifurcation.

**Main Theorem:** *There exists a codimension-one bifurcation of vector fields $X_\lambda$ having the property that for $\lambda < 0$, $X_\lambda$ is a Morse-Smale vector field, and for $\lambda > 0$, $X_\lambda$ is a Smale vector field whose closed orbits span every possible knot and link type.*

The starting point for this bifurcation is a partially hyperbolic saddle-node bifurcation, originally described by Shil'nikov [27].[1] In this bifurcation, a bouquet of homoclinic curves to a saddle-node fixed point is fattened up to a nontrivial hyperbolic invariant set. This bifurcation was significant in that it showed that a multi-dimensional system might change from simple to complex behavior via a codimension one bifurcation and it has been of interest in applications because it is believed to play a role in intemittency [24]. An explicit system of ODE's exhibiting the Shil'nikov bifurcation was produced in [8].

Recall that a *Morse-Smale* vector field is a vector field whose chain-recurrent set consists solely of a finite collection of hyperbolic fixed points and closed orbits; hence, in such a field, the *link* of closed orbits is particularly simple. On the other hand, a *Smale* vector field is one whose chain-recurrent set breaks into a finite collection of fixed points, closed orbits, and nontrivial saddle sets, all of which are hyperbolic. In a Smale flow which is not Morse-Smale, there are an infinite number of distinct closed orbits in each saddle set, and the associated link is topologically complex [10].

The question of existence of Smale flows having *all* knots and links as closed orbits was originally raised by M. Hirsch [30] and later popularised in a conjecture of Birman and Williams [3]. The issue was resolved in [12], where several examples were produced using the theory of templates. In [13], a 1-parameter family of ODEs leading to all knots and links was described; however, the bifurcation was by no means instantaneous, but rather, the knot types built up gradually. The results of this paper show that such a continuous gradation is unnecessary. Several authors have considered the knot-theoretic data associated to Morse-Smale flows [21, 29, 7] and bifurcations thereof [6].

For definitions of dynamical terms, see any of the excellent modern texts, including [15, 25]. Techniques in global bifurcation theory used in §2 can be found in [1, 23]. Comprehensive introductions to knot theory and template theory are, respectively, [26] and [14].

---

[1] Following a seminar on the Shil'nikov bifurcation given by T.Y. at Georgia Tech in 1992, K. Mischiakow proposed that this bifurcation might be analyzed for the knot structure of its periodic orbits.

# 1 Introduction

Recall that a *knot* is an embedded circle $K \subset S^3$, where $S^3$ denotes the 3-sphere[2]. Two knots are said to be *isotopic* if there exists a continuous one-parameter family of knots deforming one to the other. An *unknot* is any knot which is isotopic to an embedded curve on $S^2 \subset S^3$. A *link* is a finite collection of disjoint knots.

The observation which allows one to apply knot-theoretic perspectives to dynamical systems is that, for a flow on $S^3$, periodic orbits are knots. When the flows under consideration are dynamically complex, the knot and link types of the associated periodic orbits are likewise difficult to understand [10]. However, a construction of R. F. Williams greatly helps in the case of a nontrivial hyperbolic saddle set.

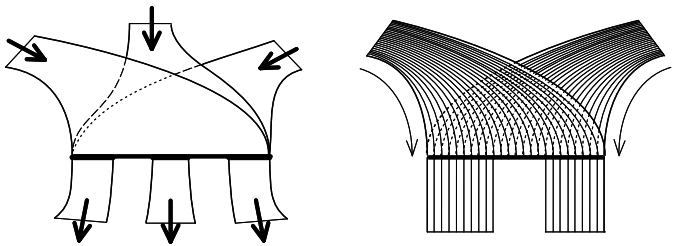

Figure 1: A branch line chart with three strips (left); each chart is fitted with an overflowing semiflow (right).

**Definition 1.1** A *template* is a compact branched 2-manifold with boundary, fitted with a smooth, expansive semiflow.

In other words, a template is what one gets by gluing together a finite number of *branch line charts* end to end, respecting the flow-direction, as indicated in Figure 1. The simplest example of a template is the *Lorenz template*, which has one such 2-strip chart. This template was introduced and analysed in [2] as a model of the attractor for the Lorenz system [20]. This template is illustrated in Figure 2(a). Also appearing is a picture of the *horseshoe template* [18]: this template likewise has one 2-strip chart, but the gluing map on the right strip is orientation-reversing. Part (c) of the figure displays a template with two branch lines.

Any template has on it a semiflow (a one-way flow) which is "overflowing" on the gaps in the branch lines; nevertheless, there are always a countable collection of disjoint periodic orbits (*cf.* the horseshoe map). We illustrate some examples in Figure 2(b). The reason why the templates of Figure 2 are named after well-known dynamical objects is explained by the Template Theorem of Birman and Williams:

**Theorem 1.2 (Birman and Williams [3])** *Given a flow on a 3-manifold $M$ with a nontrivial hyperbolic invariant set $\Lambda$, the periodic orbits of $\Lambda$ are in bijective correspondence with the set of closed orbits on some embedded template $\mathcal{T} \subset M$. The correspondence preserves all knot and link types.*

The idea is that one collapses out a foliation consisting of 1-d stable manifolds in order to compress a three-dimensional neighborhood of $\Lambda$ down to a branched surface in $M$. The identification of points with identical asymptotic futures does not alter periodic orbits.

---

[2] The reader unaccustomed to working in $S^3$ may think of it as $\mathbb{R}^3$ without loss with respect to this paper.

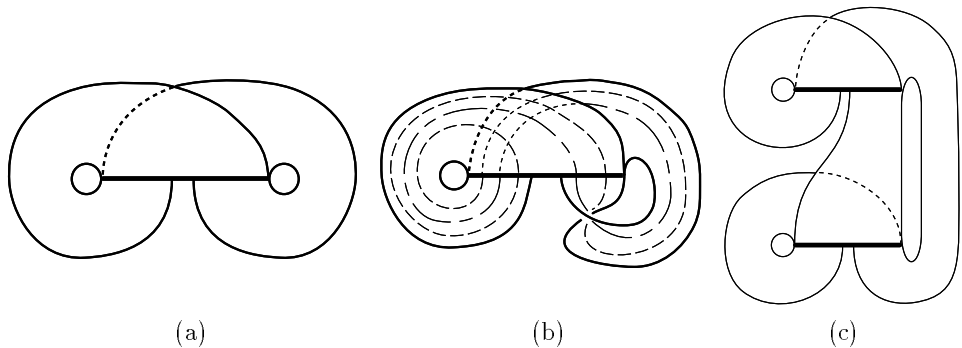

Figure 2: (a) The Lorenz template; (b) The horseshoe template with a pair of closed orbits; (c) The universal template $\mathcal{V}$.

The template theorem allows one to approximate saddle sets and attractors with an ostensibly simple combinatorial object, which may then be analysed to understand knotting and linking of closed orbits. This procedure has been carried out successfully in several model systems: see [2] for the Lorenz system, and note in particular the work of Holmes [18, 16, 17]. In all cases considered, the knots and links on a given template, though numerous, were always some restricted class of knot types.

In [11, 12], it was shown that there exist templates embedded in $\mathbb{R}^3$ in such a way as to contain *every* type of knot and link: that is, given any arbitrary link in $\mathbb{R}^3$, there is an isotopic copy of it on the template as closed orbits.

**Definition 1.3** A *universal template* is an embedded template $\mathcal{T} \subset \mathbb{R}^3$ among whose closed orbits can be found representatives of every knot and link isotopy class.

**Theorem 1.4 (Ghrist [11])** *The template $\mathcal{V}$ pictured in Figure 2(c) is universal.*

## 2 The Shil'nikov Saddle-node

In this section, we describe a class of one-parameter families which exhibit the Shil'nikov saddle-node bifurcation by introducing a set of local and global conditions. Denote by $X_\lambda$, $\lambda \in [-\lambda_0, \lambda_0]$ a one-parameter family of vector fields on a 3-dimensional manifold $M^3$ having flow $\Phi_\lambda^t$. We first assume that:

**X1** For each $\lambda$, $X_\lambda$ is $C^k$, $k \geq 1$ and $X_\lambda$ depends smoothly on $\lambda$.

### 2.1 Local saddle-node arcs

We next specify that $X_\lambda$ satisfy certain local conditions.

**X2** $X_0$ has a saddle-node singular point $p$ which is a saddle in the hyperbolic directions.

**X3** $X_\lambda$ generically unfolds the saddle-node singularity.

By **X2** we mean that $X_0(p) = 0$ and the spectrum of $DX_0(p)$ is equal to $\{\alpha, 0, -\beta\}$ where $\alpha, \beta > 0$. Associated with these eigenvalues are the eigenspaces in $T_p M$ denoted $E_p^u$, $E_p^c$, and $E_p^s$ — the unstable, center, and stable subspaces respectively. Within a sufficiently small neighborhood of $p$ there exists a unique forward- (backward-) invariant manifold $W_{lc}^s(p)$ ($W_{lc}^u(p)$) such that $T_p W^s_{lc}(p) = E^s_p$ ($T_p W^u_{lc}(p) = E^u_p$). This is commonly known as the local stable (local unstable) manifold. It is well known that there are also (non-unique) local invariant center manifolds tangent to $E^c_p$. By a saddle-node in **X2** we mean that the flow restricted to any of the center manifolds has a quadratic degeneracy at $p$. (In [4] it was shown that the flows on any two center manifolds are smoothly conjugate.)

Any sufficiently small $C^k$ neighborhood $V$ of $X_0$ contains a codimension-one hypersurface $V^{sn}$ where $X \in V^{sn}$ also has a saddle-node singular point in $U$. The set $V^{sn}$ divides $V$ into two disjoint regions $V^-$ and $V^+$, where each $X \in V^-$ has a pair of hyperbolic singular points in $U$ and every $X \in V^+$ has no singular points in $U$. By **X3** we mean that the arc $\{X_\lambda, -\lambda_0 < \lambda < \lambda_0\}$ is transversal to $V^{sn}$ at $X_0$.

Equivalent to **X2** and **X3**, we may assume that there exist local coordinates $(x, y, \theta)$ in $U$ where $p$ is given by $(0, 0, 0)$, for which $X_\lambda|_U$ has the form (possibly after reparameterization):

$$\begin{aligned}
\dot{x} &= (\alpha + f(x, y, \theta, \lambda))x \\
\dot{y} &= (-\beta + g(x, y, \theta, \lambda))y \\
\dot{\theta} &= \lambda + \theta^2 + h(x, y, \theta, \lambda),
\end{aligned} \tag{1}$$

where at the origin, $f = g = h = h_\theta = h_{\theta\theta} = 0$ (For normal forms see [28] and [19]).

For $\lambda < 0$ the vector field has two hyperbolic singular points, $p_\lambda$ and $q_\lambda$, in $U$. As $\lambda \nearrow 0$, these points coalesce at $p$. With the convention that the $\theta$-coordinate of $q_\lambda$ is greater than that of $p_\lambda$, the Stable Manifold Theorem implies that $p_\lambda$ has a 2-dimensional local stable manifold $W^s_{lc}(p_\lambda)$ and a 1-dimensional local unstable manifold $W^u_{lc}(p_\lambda)$. Likewise, the local stable manifold of $q_\lambda$ is 1-dimensional and the local unstable manifold has dimension 2. Standard results on local stable and unstable manifolds (for example [9]) yield that $W^u_{lc}(q_\lambda)$ and $W^s_{lc}(p_\lambda)$ intersect transversally along a unique local heteroclinic solution curve $\gamma_{lc,\lambda}$. In addition, one has $W^s_{lc}(q_\lambda) = \partial \overline{W^s_{lc}(p_\lambda)} \cap U$, and $W^u_{lc}(p_\lambda) = \partial \overline{W^u_{lc}(q_\lambda)} \cap U$, as in Figure 3 (left).

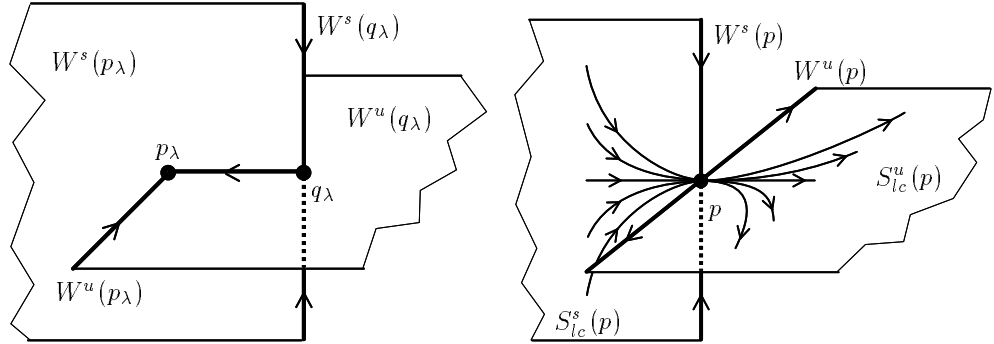

Figure 3: (left) The local structure of stable and unstable invariant manifolds for $\lambda < 0$, and (right) for $\lambda = 0$.

As mentioned previously, for $\lambda = 0$ we have $\dim(W^s_{lc}(p)) = \dim(W^u_{lc}(p)) = 1$. In addition, **X2** implies that there are a pair of local invariant manifolds (with boundary) associated to $p$, denoted $S^s_{lc}(p)$ and $S^u_{lc}(p)$, the local stable and unstable sets respectively[3]. These are defined by:

$$S^s_{lc}(p) = \{X \in U : \Phi^t(X) \in U, \forall t > 0, \text{ and } \Phi^t(X) \to p, \text{ as } t \to +\infty\},$$

---

[3] We note that some authors use $W^s_{lc}$ to denote the stable set (and call it the stable manifold) and $W^{ss}_{lc}$ to denote the stable manifold (and refer to it as the *strong stable* manifold). The notation we choose to follow comes from [1].

$$S^u_{lc}(p) = \{X \in U : \Phi^t(X) \in U, \forall t < 0, \text{ and } \Phi^t(X) \to p, \text{ as } t \to -\infty\}.$$

The relationship between the (un)stable manifolds and the (un)stable sets is that (see [22]) $S^s_{lc}(p) = \lim_{\lambda \nearrow 0} W^s_{lc}(p_\lambda)$, and $S^u_{lc}(p) = \lim_{\lambda \nearrow 0} W^u_{lc}(q_\lambda)$, where the limit denotes the topological limit. Also, the boundary of $S^s_{lc}(p)$ is $W^s_{lc}(p)$ and $\partial S^u_{lc}(p) = W^u_{lc}(p)$.

The local stable set is obviously unique and is necessarily included in any 'local center-stable manifold', *i.e.*, in any invariant manifold tangent to $E^c_p \oplus E^s_p$. Because of the exponential contraction in the $E^s$ direction, the orbits of $X_0$ on $S^s_{lc}(p) \setminus W^s_{lc}(p)$ all approach $p$ along the $E^c_p$-direction. The same holds for backward orbits on $S^u_{lc}(p) \setminus W^u_{lc}(p)$, as illustrated in Figure 3 (right).

### 2.2 Global conditions

Extend $S^s_{lc}(p)$, $S^u_{lc}(p)$, $W^s_{lc}(p)$, and $W^u_{lc}(p)$ by the flow to global invariant sets $S^s(p)$, $S^u(p)$, $W^s(p)$, and $W^u(p)$, respectively. That is, let $S^s(p) = \{x | \Phi^t(x) \in S^s_{lc}(p), \text{ for some } -\infty < t < \infty\}$, and similarly for $S^u(p)$, $W^s(p)$, and $W^u(p)$.

The next condition was first considered by Shil'nikov [27].

**X4** $S^u(p)$ intersects $S^s(p)$ transversally along a finite collection of orbits $\{\gamma_i\}|^n_{i=1}$, none of which are in $W^u(p)$ or $W^s(p)$.

Denote by $\Gamma_0$ the bouquet of homoclinic loops $\Gamma_0 = \{p\} \cup \{\bigcup^n_{i=1} \gamma_i\}$.

**Theorem 2.1 (Shil'nikov [27])** *Given conditions* **X1 - X4** *and a small fixed neighborhood $U$ of $\Gamma$, there exists a $\overline{\lambda} > 0$ such that for all $0 < \lambda < \overline{\lambda}$, the flow $\Phi_\lambda$ restricted to the invariant set of $U$ is topologically conjugate to a suspension of a topological Bernoulli shift on $n$ symbols.*

For each $\lambda$, let $\Gamma_\lambda$ denote the invariant set of $U$. By transversality and the genericity of the unfolding conditions, $\Gamma_\lambda = \{p_\lambda\} \cup \{q_\lambda\} \cup \{\bigcup^n_{i=1} \gamma_{i,\lambda}\}$ for $\lambda < 0$. It is shown in [31] and [19] that for $\lambda > 0$, the invariant set $\Gamma_\lambda$ is hyperbolic with respect to the flow. This hyperbolicity along with Theorem 1.2 implies that the periodic orbits of $\Gamma_\lambda$ correspond to those on an embedded template. The precise embedding of the template is determined by the structure of invariant sets to $\Gamma_0$.

### 2.3 Pretemplates

The transversality condition of **X4** that $\Gamma_0$ is a normally hyperbolic set and hence its local unstable set $S^u_{lc}(\Gamma_0)$ is the union of local unstable manifolds of its points. This set is of particular importance to the topology of the flow and will hence be called the *pretemplate*.

**Proposition 2.2** *The pretemplate $S^u_{lc}(\Gamma_0)$ is a compact branched two-manifold with a unique branch line along which the flow is invariant.*

*Proof:* By definition of $\Gamma_0$, the pretemplate equals

$$S^u_{lc}(\Gamma_0) = W^u_{lc}(p) \cup \left(\bigcup^n_{i=1} W^u_{lc}(\gamma_i)\right).$$

Outside a neighborhood of $p$, each $\gamma_i$ is a bounded distance away from every other $\gamma_j$; hence, $S^u_{lc}(\Gamma_0)$ is a disjoint set of $n$ compact strips outside a neighborhood of $p$. Transversality and the Lambda Lemma [23] imply that every $W^u_{lc}(\gamma_i)$ is tangent to the one-dimensional unstable manifold $W^u_{lc}(p)$ at $p$. However, in the limit as $t \to +\infty$, every $\gamma_i(t)$ is contained within the local stable manifold $W^s_{lc}(p)$, which is transversal to the unstable eigendirection. Hence, in the forwards limit, the strips $W^u_{lc}(\gamma_i)$ are mutually tangent and "stacked" onto $W^u_{lc}(p)$, which forms the (unique) branch line for the pretemplate. On the other hand, in the limit as $t \to -\infty$, $W^u_{lc}(\gamma_i)$ limits onto $W^u_{lc}(p)$, but the uniqueness of $S^u(p)$ implies that $W^u_{lc}(\gamma_i)$ coincides with $S^u_{lc}(p)$ near $p$ and thus separates as the flow leaves a neighborhood of $p$. Thus we obtain the $n$-strip branch line chart. □

Note that along the branch line $W^u_{lc}(p)$, the flow is invariant and in fact is expanding, split by the fixed point $p$; hence, the pretemplate is not a template in the sense of Definition 1.1. However, it is a branched two-manifold and thus diffeomorphic to a template.

**Proposition 2.3** *For $\lambda > 0$ sufficiently small, the template obtained by collapsing the local stable manifold of orbits in $\Gamma_\lambda$ is isotopic as a branched manifold to the pretemplate $W^u_{lc}(\Gamma_0)$.*

*Proof:* Let $U$ be a neighborhood of $\Gamma_0$ which consists of an $\epsilon$-ball $B$ centered at $p$ and $\delta$ neighborhoods $\{G_i\}_{i=1}^n$ of the homoclinics where $\delta$ is sufficiently small that these neighborhoods are disjoint outside of $B$. In [31] it was shown that the stable manifolds of $\Gamma_\lambda$ are locally $C^k$ close to the stable manifolds of $\Gamma_0$. Thus on each $G_i \setminus B$ the template associated with $\Gamma_\lambda$ obtained by collapsing the stable direction is simply a strip which is a smooth deformation of $W^u_{lc}(\Gamma_0) \cap G_i$. The proof of Theorem 2.1 implies that, within $B$, the part of $\Gamma_\lambda$ entering from each $G_i$ is stretched by the hyperbolicity across every $G_j$. Hence, each strip in the template must connect to every other strip inside $B$. Monotonicity of the flow in $B$ gives that this can only be accomplished via a local chart equivalent to a single branch line with $n$ strips. Thus the template is a smooth perturbation of the pretemplate. □

## 3 Twist and Pleatings

The knot-theoretic results we derive in the following sections depend crucially on the notion of *twist*. For simplicity in the definitions, we restrict our systems to those in $\mathbb{R}^3$ (or, equivalently, $S^3$).

**Definition 3.1** Let $\gamma$ be a simple closed curve on an annulus or Möbius strip $A$ embedded in $\mathbb{R}^3$. The *push-off* of $\gamma$ with respect to $A$ is the embedding of $\nu^1(\gamma)$, the unit normal bundle to $A$ restricted to $\gamma$, in the complement $\mathbb{R}^3 \setminus \gamma$. The *twist* of $\gamma \subset A$ is then defined as the homology class of the push-off in $H_1(\mathbb{R}^3 \setminus \gamma; \mathbb{Z}) \cong \mathbb{Z}$. Equivalently, it is the linking number of $\gamma$ with its pushoff(s) [26].

Note that if $A$ is an annulus, the push-off consists of two disjoint circles (one on each side), whereas if $A$ is a Möbius band, then the push-off is connected, but wraps about $A$ twice. One must establish a sign convention for computing twists — this is equivalent to choosing a generator for $H_1(\mathbb{R}^3 \setminus \gamma; \mathbb{Z}) \cong \mathbb{Z}$. We will work under the convention that *left-over-right is positive*, so that the twist associated to the core of the twisted strip of the horseshoe template of Figure 2(b) is $+1$.

From the definition, it follows that the twist of a closed orbit $\gamma$ on a template $\mathcal{T}$ is the twist associated to a "ribbon" transverse to $\mathcal{T}$ along $\gamma$. Likewise, the twist associated to a homoclinic orbit with well-defined stable and unstable manifolds (such as those considered in §2) is the twist with respect to the local stable or unstable manifolds.

As in the construction of §2, let $\Gamma_0$ denote the set of homoclinic orbits to $p$ along with $p$. Upon orienting the unstable set $S^u(p)$, there is a natural ordering on the homoclinic orbits $\gamma_i$ as follows: intersect $S^u_{lc}(p)$ with a small 2-sphere centered at $p$ to obtain an oriented arc transverse to $\Gamma_0$. For the remainder of this paper, we will assume that the set $\{\gamma_i\}$ is given this fixed ordering.

**Definition 3.2** The *twist signature* of the ordered homoclinic curves $\{\gamma_i\}_1^n$ is defined as the ordered sequence of integers $\tau = (\tau_1, \tau_2, \ldots, \tau_n)$, where $\tau_i \in \mathbb{Z}$ is the twist of $\gamma_i$.

Since we order the $\gamma_i$ curves following the unstable set $S^u(p)$, the twist signature is related to how $S^u(p)$ wraps around itself, as seen in a cross-section. A natural class of systems to consider are those which have a simple cross-sectional behavior.

**Definition 3.3** In the construction of §2, we say that the system is *pleated* if there exists a smooth embedded disc $D \subset M$ such that (1) $D$ is transverse to $S^u(p)$ and $S^s(p)$; (2) $D$ intersects every component $\gamma_i$ of $\Gamma_0$ once; and (3) $D \cap S^s(p)$ and $D \cap S^u(p)$ are each a connected, proper submanifold of $D$.

The essence of this definition is expressed in Figure 7 (left), where it is seen that in the section induced by $D$, the unstable set $S^u(p)$ crosses the stable set $S^s(p)$ in an alternating fashion, forcing it to fold or "pleat" nicely.

**Definition 3.4** The *pleating signature* of a pleated system is defined to be the isotopy class of the graph $D \cap \{S^u(p) \cup S^s(p) \cup \Gamma_0\}$ within $D$.

**Lemma 3.5** *The pleating signature is well defined.*

*Proof:* Given a disc $D$ defining a pleating, one may deformation retract $D$ rel $\{S^s(p) \cup S^u(p) \cup \Gamma_0\}$ so that $\partial D$ is close to $S^s(p) \cup S^u(p)$. Then, shrink the four "ends" of $D$ where the branches of $S^u(p)$ and $S^s(p)$ intersect $\partial D$, until each one is very close to $\Gamma_0$. This reduces the pleating disc to a "normal form" without affecting the pleating signature. Given any other disc $D'$, reduce it likewise and use the action of the flow to isotope $D$ to $D'$ as follows.

There is a natural cellular decomposition on $D$ and $D'$ induced by the invariant manifolds, where the 0-skeleton is given by intersection with $\Gamma_0$, and the 1-skeleton is given by intersection with $S^u(p)$ and $S^s(p)$. Since $D$ and $D'$ intersect each $\gamma_i$ once, there is a unique direction in which to flow the 0-skeleton of $D$ to that of $D'$ along $\Gamma_0$ without passing through $p$. Each connected component of $S^u(p)$ and $S^s(p)$ is an arc connecting distinct points on the 0-skeleton. Using the flow directions induced on $\Gamma_0$, extend the flow on the 0-skeleton of $D$ to the 1-skeleton. There is a natural extension over the 2-skeleton, completing the isotopy from the pleating signature of $D$ to that of $D'$. □

**Corollary 3.6** *All of the homoclinic curves $\gamma_i$ in a pleated system are isotopic as knots.*

*Proof:* The flow-isotopy of Lemma 3.5 implies that there is a 1-parameter family of discs, each transversal to all the $\gamma_i$ with one intersection. This family extends only up to a neighborhood of $p$, where, however, the local behaviour at $p$ is completely determined: each $\gamma_i$ approaches $p$ as in Figure 3 (right). Thus, filling in with the local picture, there is a (non-invariant) solid torus $V \subset \Gamma_0$ such that each $\gamma_i \cup \{p\}$ is isotopic to the core longitudinal curve, which defines the knot type of all the $\gamma_i$ curves. □

**Proposition 3.7** *Given a pleated system in $\mathbb{R}^3$ as above, the isotopy class of the associated pretemplate is determined by the pleating signature, the twist signature, and the knot type of the homoclinic curves $\{\gamma_i\}$.*

*Proof:* From §2.3, the pretemplate $S^u_{lc}(\Gamma_0)$ consists of a single branch line and $n$ strips, one along each $\gamma_i$. Fixing the knot types of the $\gamma_i$ establishes a rough equivalence between individual strips. Two other pieces of data determine the isotopy type of the template. The first is the twist in each strip: this is determined by the twist signature. The second is the relative order with which the strips emanate from and overlap into the branch line: this is determined by the pleating signature. □

One reason for restricting to the class of pleated systems is that their twist signatures behave nicely.

**Lemma 3.8** *The twist signature of a pleated unstable set is* incremental: *that is, $\tau_{i+1} = \tau_i \pm 1$ for all $i$.*

*Proof:* Fix some $i$ and consider the homoclinic curves $\gamma_i$ and $\gamma_{i+1}$, having twist $\tau_i$ and $\tau_{i+1}$. The pushoff with respect to $S^u(p)$ is a subset of $S^s_{lc}(\gamma_i)$. Hence, to compute the twist of $\gamma_{i+1}$, note that in the limit as $t \to -\infty$, the pushoffs to $\gamma_i$ and $\gamma_{i+1}$ coincide within $S^s(p)$. The pleating signature indicates the folding of $S^u(p)$ as one sends $t \to \infty$: there is a single "U" fold by the Jordan Curve Theorem. Hence, in moving from $\gamma_i$ to $\gamma_{i+1}$ by sliding along $D \cap S^u(p)$, one gains an extra single twist in the pushoff. □

## 4 Universal Templates

Several properties of universal templates are derived in [12] and [14]. However, the question of when a template is universal is quite subtle and has not yet been fully answered. For example, the horseshoe template of Figure 2(b) does not even support all torus knots (see [3, 18]), yet by reversing the twist in the right strip, the template becomes universal. In fact, *every* template can be re-embedded so as to be universal [14]. The fundamental difference lies in the *twist* of the strips in the template. The technical theorem that we use to classify which of the bifurcations of §2 lead to universal templates is the following:

**Theorem 4.1 (Ghrist [12],[14](p. 91))** *Let $\mathcal{T}$ denote an embedded template in $\mathbb{R}^3$. Sufficient conditions for $\mathcal{T}$ to be universal are as follows. There exists a pair of disjoint closed orbits on $\mathcal{T}$, $\kappa$ and $\kappa'$, such that: (1) they are separable (unlinked) unknots; (2) $\tau(\kappa) = 0$ and $\tau(\kappa') \neq 0$, where $\tau$ denotes twist; and (3) There exists a subset of $\mathcal{T}$ isotopic to Figure 4(a), in the case where $\tau(\kappa') < 0$; or, isotopic to Figure 4(b), in the case $\tau(\kappa') > 0$.*

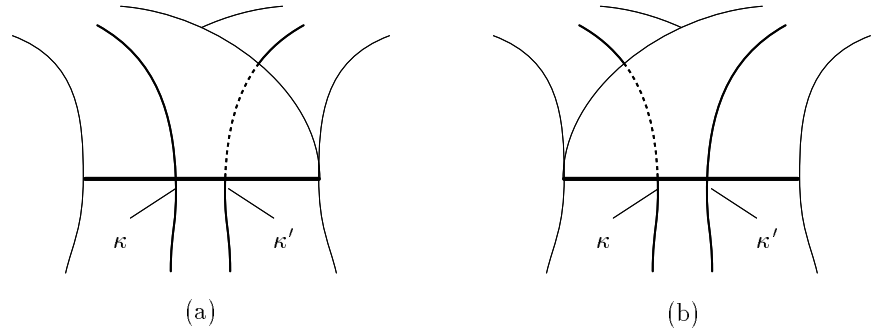

Figure 4: (a) A positive branch line crossing; (b) A negative branch line crossing.

To derive sufficient conditions for a pleated system to give rise to a universal template, the way in which the flow is embedded is crucial. Two features of the embedding are crucial: whether the individual $\gamma_i$ are knotted, and whether they are mutually linked via twisting. The restrictions due to the former are very easy to understand, whereas the restrictions due to twisting are more subtle.

**Lemma 4.2** *Given a pleated system, if the homoclinic curves $\gamma_i$ are knotted, then the associated template for $\lambda > 0$ does not support any unknotted closed orbits, and hence is not universal.*

*Proof:* The *genus* of a knot is a classical invariant for knots defined to be the minimal genus (number of handles) for an orientable surface whose boundary is the knot [26]. The proof follows easily using the formula in [5, p. 20] to compute that the genus of every orbit on the template is bounded below by the genus of each $\gamma_i$. □

**Theorem 4.3** *Let $X_\lambda$ be a family of vector fields unfolding a pleated Shil'nikov saddle-node in $\mathbb{R}^3$. Then the pretemplate at $\lambda = 0$ is isotopic to a universal template for $\lambda > 0$ if (1) the homoclinic orbits $\{\gamma_i\}$ are unknots; and (2) the twist signature of the system has both positive and negative terms.*

*Proof:* If the twist signature has mixed signs, Lemma 3.8 implies that there is a portion of the pleating where the twist signature looks like either $(\ldots, -1, 0, 1, \ldots)$ or $(\ldots, 1, 0, -1, \ldots)$. We may restrict the associated template to that generated by these three strips alone, since any other strip only *adds* closed orbits to the system: denote by $\mathcal{S}$ this subtemplate.

By Proposition 3.7, $\mathcal{S}$ is completely determined up to isotopy by the twist signature, the pleating signature, and the knot types of the homoclinic curves. The twist signature is either $(-1, 0, 1)$ or $(1, 0, -1)$ and the pleating signature is likewise restricted to be one of four possible cases as illustrated in Figure 5, only two of which are compatible with each twist signature. By hypothesis, we have specified that the $\gamma_i$ curves are unknotted. Hence, there are precisely four different possibilities for the subtemplate, up to isotopy. However, there is a symmetry in the pleatings which corresponds to taking a *mirror image* of the template (reversing the sign of all crossings). Since a template is universal if and only if its mirror image is universal, we can reduce to two cases.

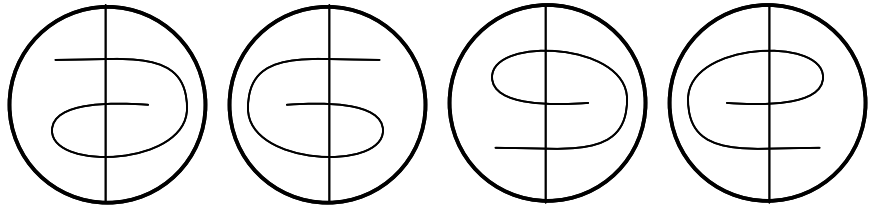

Figure 5: The possible pleatings associated to the twist signatures $(\ldots, -1, 0, 1, \ldots)$ and $(\ldots, 1, 0, -1, \ldots)$.

We use the sufficient conditions of Theorem 4.1 to show that, in each case, the template is universal. In Figure 6, we illustrate the template associated to each pleating type (here, the top and the bottom of each diagram is identified), along with a pair of curves $\{\kappa, \kappa'\}$ in each. It is straightforward to verify that these curves are separable unknots with twist 0 and 6 respectively meeting at a negative branch line crossing as per Figure 4(b). By Theorem 4.1, these conditions guarantee that the subtemplate $\mathcal{S}$ is universal. □

**Remark 4.4** It is not true that the presence of mixed signs in the twist signature is necessary for the existence of a universal template. Consider Figure 7 (left), in which is shown an oriented pleating diagram. By setting the twist of the first $\gamma_i$ to be zero, this system has twist signature $\tau = (0, 1, 0, 1, 2, 1, 0, 1)$, which is nonnegative. Taking the associated template and ignoring all strips with nonzero twist, we obtain the subtemplate of Figure 7 (identify top and bottom). As shown, there exists a pair of separable unknots with twist 0 and $-2$; this is a universal template as per Theorem 4.1.

The reason for the presence of a universal template is, as it must be, the existence of mixed positive and negative twisting in the pleating. While Theorem 4.3 guarantees that mixed twisting in the strips is sufficient, we see that the twisting may be hidden within the pleating itself. It remains an open problem to find necessary and sufficient conditions for a pleated system to support a universal template.

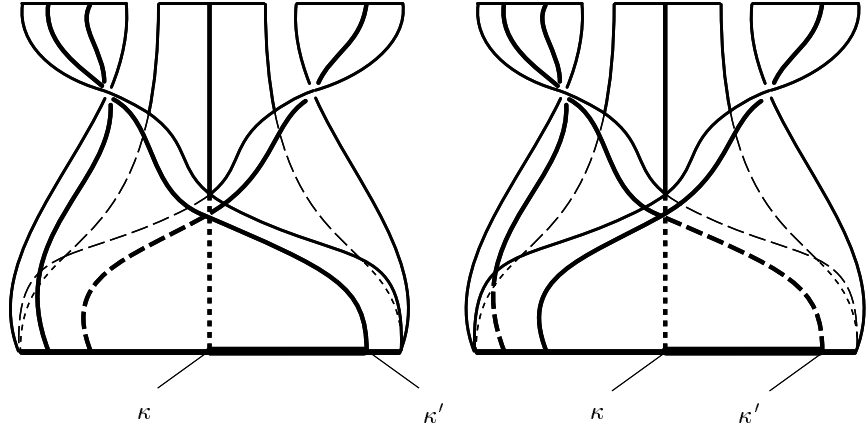

Figure 6: The subtemplates $\mathcal{S}$ associated to the twist signatures $(1, 0, -1)$ and $(-1, 0, 1)$ [identify tops and bottoms], along with pairs of separable unknots forcing all other knots.

## 5  From Morse-Smale to a Universal Template

In the remainder of this work, we construct an explicit example of a pleated system which satisfies the template-theoretic conditions of §4.

Let $\Psi^{sn}$ denote the codimension one hypersurface of saddle-node vector fields in the space of smooth vector fields. A 'versal deformation' of a $k$-degenerate vector field $\tilde{X}$ is a $k$-parameter family of vector fields which is transversal at $\tilde{X}$ to the codimension $k$ bifurcation hypersurface containing $\tilde{X}$.

**Theorem 5.1** *There exists an open set $V \subset \Psi^{sn}$ such that, any $\tilde{X} \in V$ has a codimension one degeneracy and if $\tilde{X}_\lambda$ is any versal deformation of $\tilde{X} \in V$, then the induced flow $\tilde{\Phi}_\lambda$ is M-S for $\lambda < 0$ and is Smale containing a universal link for $\lambda > 0$.*

*Proof:* We begin by prescribing the flow on a solid torus $M \cong D^2 \times S^1$. It is a simple exercise to embed the flow into $\mathbb{R}^3$ and $S^3$, preserving Morse-Smale. A suggestive illustration is given as Figure 8: we build up this picture via additional assumptions.

**X5** On $M \setminus \{p\}$, the flow $\Phi_0$ is monotonic in the longitudinal direction.

That is, there are coordinates on the solid torus $(x, y, \theta)$ such that $\dot{\theta} > 0$ except at $p$. (We will take $p$ to have local coordinates $(0, 0, 0)$.)

One of the consequences of **X5** is of a local nature. Namely, let $\Sigma^-$ and $\Sigma^+$ be cross-sections of the torus close to $p$ given by $\{\theta = -\overline{\theta}\}$ and $\{\theta = \overline{\theta}\}$, where $\overline{\theta} > 0$ is small. Let $U$ be the neighborhood of $p$ bounded by $\Sigma^-$ and $\Sigma^+$. Condition **X5** implies that the two branches of $W^u_{lc}(p)$ are bent forward, so that they intersect $\Sigma^+$ in two points and $\overline{S^u_{lc}(p)}$ intersects $\Sigma^+$ in an arc connecting these points. The endpoints of the arc $a$ and $b$ are precisely the intersection of $\overline{W^u_{lc}(p)}$ with $\Sigma^+$. Similarly, $\overline{S^s_{lc}(p)} \cap \Sigma^-$ is an arc with endpoints we denote $c$ and $d$. In the local coordinates for which $X_0$ is given by (1), each of these arcs has the form of a line segment. We denote these arcs by $\overline{ab}$ and $\overline{cd}$, as illustrated in Figure 8.

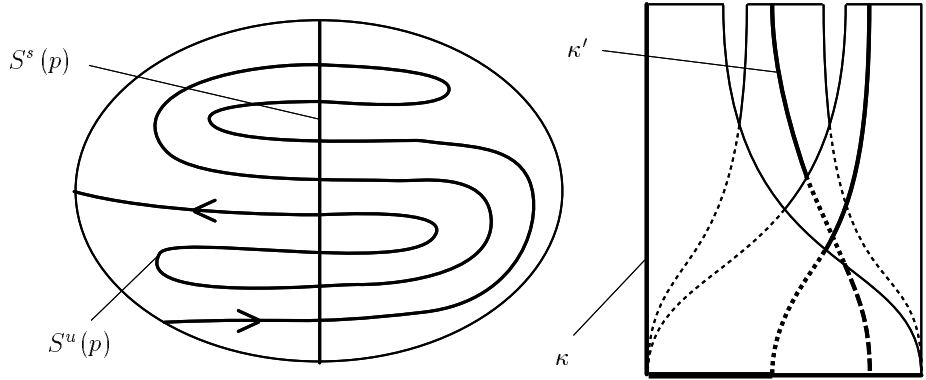

Figure 7: A certain pleated system (left) yields a subtemplate (right) containing a pair of unknots which forces all others.

Condition **X5** may be realized explicitly in $U$ by a vector field of the form in Equation 1, with $f = g = 0$ and $h = \mu x^2 + \nu y^2$, where $\alpha$, $\beta$, $\mu$, and $\nu$ are all positive. The Stable Manifold Theorem implies that $W_{lc}^u(p)$ is given by $(x, 0, H(x))$ where $H(x)$ is a smooth function with $H(0) = H'(0) = 0$. Substitution of $\theta = H(x)$ into these equations gives: $\alpha x H'(x) = H(x)^2 + \mu x^2$; hence, $\alpha x H'(x) \geq \mu x^2$. This, along with $H(0) = H'(0) = 0$, implies $H(x) > \mu x^2/2\alpha$. Thus, by adjusting the ratios $\mu/\alpha$ and $\nu/\beta$ we may 'bend' $W_{lc}^u(p)$ and $W_{lc}^s(p)$ by arbitrary amounts.

By the above observations, $S^u(p)$ has the shape of a 'ribbon' as it leaves $U$, the two branches of $W^u(p)$ forming its edges. Note that by the continuous dependence of $\Phi_\lambda$ on $\lambda$, $W^u(p_\lambda)$ is also bent in $U$, so that $W^u(q_\lambda)$ is also shaped like a ribbon, as in Figure 8.

Let $\Sigma^0$ be the coordinate cross section in $U$ which contains $p$, that is, $\Sigma_0 = \{\theta = 0\}$. The flow induces various mappings between the sections $\Sigma^-$, $\Sigma^+$, and $\Sigma^0$. Namely, let $F_{gl} : \Sigma^+ \to \Sigma^-$ be the diffeomorphism induced by $\Phi_0^t|_{M \setminus U}$. That is, for any point $x \in \Sigma^+$ let $F_{gl}(x)$ be the first (as ordered by $t$) intersection of $\Phi^t(x), t > 0$ with $\Sigma^-$. We define a map $F_- : \Sigma^- \to \Sigma^0$ in a similar fashion, except for $x \in \overline{cd}$. There, since $\Phi^t(x) \to p$ as $t \to +\infty$ we take $F_-(x) = p$. Thus defined, $F_-$ is differentiable, though not invertible since $F_-(\overline{cd}) = p$. Finally, we define the map $F_+ : \Sigma^+ \to \Sigma^0$ induced by $\Phi^{-t}$, where it is seen that $F_+(\overline{ab}) = p$. The map $F_+$ is invertible on $\Sigma^+ \setminus \overline{ab}$ and henceforth when we refer to $F_+^{-1}$ we make appropriate assumptions about its domain.

Let $A$, $B$, $C$, and $D$ denote the regions of $\Sigma^-$ displayed in Figure 9 (left). Consider the image of $\overline{cd} \cup F_{gl}(\overline{ab})$ under $F_{lc} \equiv \circ F_+^{-1} \circ F_-$. The map $F_-$ collapses $\overline{cd}$ to $p$, including the points contained in $\overline{cd} \cap F_{gl}(\overline{ab})$. The action of $F_+^{-1}$ is to blow-up $p$ into the segment $\overline{ab}$. Since the boundaries of the regions $A$, $B$, $C$, and $D$ are deformed continuously, it follows that the action of $F_{lc}$ on the partition $\{A, B, C, D\}$ is as shown in Figure 9 (middle).

**X6** Let $\Phi_0|_{M \setminus U}$ be such that the action of $F_{gl} : \Sigma^+ \to \Sigma^-$ is a 'double horseshoe' map as shown in Figure 9 (right). Denote by $\gamma_1$ the homoclinic orbit adjacent to $a$ and let $\gamma_1$ have zero twist.

It follows from the structure of $\overline{cd}$ and $F_{gl}(\overline{ab})$ that the system is pleated with pleating signature derived from Figure 9 (right). By re-embedding the solid torus of the flow with meridional twisting, it is clear that the flow can be constructed so that $\gamma_1$ has zero twist without altering the pleating signature. It follows that the twist signature of the pretemplate is $\tau = (0, 1, 0, -1)$. Thus by Theorem 4.3, the invariant set for $\lambda > 0$ contains a universal link.

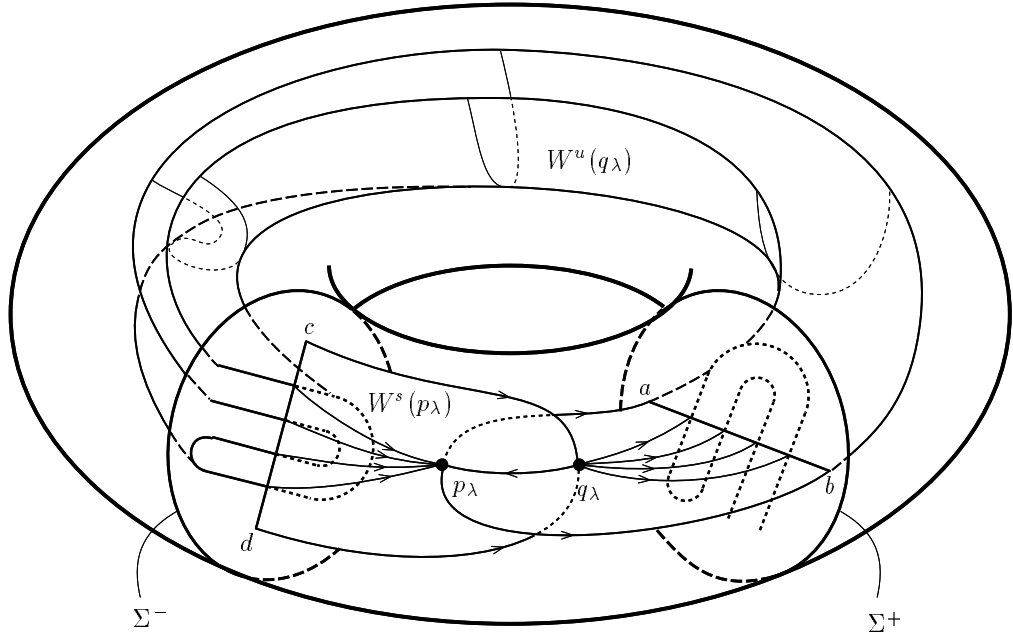

Figure 8: A pleated system which bifurcates into a universal template.

Denote $F \equiv F_{gl} \circ F_{lc}$. Note that $F(D) \subset D$ and $F^{-1}(A) \subset A$.

**X7** $F|_D$ is a contraction, as is $F^{-1}|_A$.

Condition **X7** implies that $\overline{D}$ has an attracting fixed point associated with an attracting periodic orbit $\Gamma_0^\omega$. Since $F|_{\overline{D}}$ fixes the boundary of the solid torus, this attracting orbit must lie on $\partial M$. Similarly, $\overline{A} \cap \partial M$ contains a repelling fixed point corresponding to a repelling periodic orbit $\Gamma_0^\alpha$.

**Proposition 5.2** *Under conditions* **X1** - **X7** *the flow* $\Phi_\lambda$ *is M-S on the solid torus $M$ for $\lambda < 0$.*

*Proof:* Let
$$\Gamma_\lambda = p_\lambda \cup q_\lambda \cup \gamma_{lc,\lambda} \bigcup_{j=1}^n \gamma_{j,\lambda}$$

for $\lambda < 0$. The set $\Gamma_\lambda$ is invariant and the dynamics on it are Morse-Smale. We show that the flow is globally Morse-Smale by showing that all the other dynamics are 'absorbed' in the attractor-repellor pair on the boundary of the solid torus $\partial M$. For $\lambda = 0$, since the images of $\overline{B}$ and $\overline{D}$ are contained in $\overline{D}$, $\Gamma_0^\omega$ is the $\omega$-limit set for all points in $B$ and $D$. Thus the only chain-recurrent set in either of these regions is $\Gamma_0^\omega$. If we consider the inverse map $F^{-1}$ it follows by a similar argument that the $\alpha$-limit set of all points of these regions except those on the common boundary $\overline{cd}$ must be $\Gamma_0^\alpha$. Finally the points on $\overline{cd}$ all limit to $p$ in forward time, in backward time they limit to either $p$ or $\Gamma_0^\alpha$. Thus the only chain-recurrent sets for $\lambda$ are $\Gamma_0^\omega$, $\Gamma_0^\alpha$, and $p$.

For $\lambda < 0$ we partition $\Sigma^-$ in such a way that the action of $F_\lambda$ on the partition is identical to that of $F$ on $\{A, B, C, D\}$. First we observe that the eigenvalues of the linearisation of $X_\lambda$ at $p_\lambda$ are $\{\alpha_\lambda, -\delta_\lambda, -\beta_\lambda\}$ where $\alpha_\lambda \approx \alpha$, $\beta_\lambda \approx \beta$ and $-\delta_\lambda < 0$ is approximately zero. It follows [9] that for $\lambda < 0$ sufficiently close to zero there is a unique smooth locally invariant manifold tangent to the eigenspace corresponding to $-\beta_\lambda$ at $p_\lambda$. This manifold we refer to as the strong-stable manifold of $p_\lambda$ and denote it by $W^{ss}_{lc}(p_\lambda)$. Similarly there is a unique smooth invariant manifold tangent to the strongly unstable eigen-direction at $q_\lambda$ which we denote by $W^{uu}_{lc}(q_\lambda)$. It follows from [9] that as $\lambda \nearrow 0$, $W^{ss}_{lc}(p_\lambda)$ converges smoothly to $W^s_{lc}(p)$ and $W^{uu}_{lc}(q_\lambda)$ converges smoothly to $W^u_{lc}(p)$. By the smooth convergence, $W^{ss}_{lc}(p_\lambda)$ intersects $\Sigma^-$ in two points which we call $c_\lambda$ and $d_\lambda$ which converge to $c$ and $d$, respectively as $\lambda \nearrow 0$. Also, $W^{uu}_{lc}(q_\lambda)$ intersects $\Sigma^+$ in points $a_\lambda$ and $b_\lambda$ which converge to $a$ and $b$. Let $\{A_\lambda, B_\lambda, C_\lambda, D_\lambda\}$ be a partition of $\Sigma^-$ as in Figure 9 (left), only parametrised by $\lambda$ appropriately.

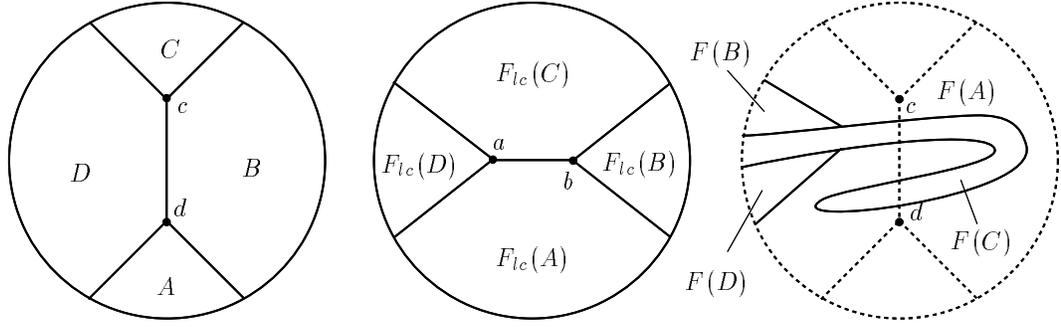

Figure 9: A partition of $\Sigma^-$ (left); The action of $F_{lc}$ on the partition (middle); and the double-horseshoe action of $F_{gl}$ (right).

The action of $\Phi|_U$ is to take $\overline{c_\lambda d_\lambda}$ to the point $p_\lambda$, $p_\lambda$ to $q_\lambda$ and $q_\lambda$ to $\overline{a_\lambda b_\lambda}$. Thus the local action on the partition is the same as that for $\lambda = 0$. The action of $\Phi_\lambda|_{M \setminus U}$ is a smooth perturbation of that for $\lambda = 0$. The transversality condition then implies that the action on the partition of $F_\lambda$, for $\lambda < 0$ sufficiently close to 0, is isotopic to that for $\lambda = 0$. It also follows that for $\lambda$ sufficiently close to zero, $F_\lambda|_D$ and $F_\lambda^{-1}|_A$ are still contractions and so there are unique attracting and repelling orbits, $\Gamma^\omega_\lambda$ and $\Gamma^\alpha_\lambda$, in $D$ and $A$ respectively. Finally note that $F_{gl}(\overline{a_\lambda b_\lambda}) \setminus \overline{c_\lambda d_\lambda} \subset B \cup D$, and so, all orbits of $S^u(p)$ which are not in $\overline{c_\lambda d_\lambda} \cap F_{gl}(\overline{a_\lambda b_\lambda})$ limit to $\Gamma^\omega_\lambda$. Thus there are no intersections of $W^u(p_\lambda)$ with $W^s(q_\lambda)$ other that $\{\gamma_{i,\lambda}\}_{i=1}^4$ and $\gamma_{lc,\lambda}$. These are transverse, $W^s(p_\lambda)$ limits to $\Gamma^\alpha_\lambda$ in backward time and $W^u(q_\lambda)$ limits to $\Gamma^\omega_\lambda$ in forward time. Thus $W^s(p_\lambda)$ and $W^u(q_\lambda)$ do not intersect.

For $\lambda < 0$, it follows from the same arguments as for $\lambda = 0$ that the chain recurrent set for $\Phi_\lambda$ consists of the periodic orbits $\Gamma^\alpha_\lambda$ and $\Gamma^\omega_\lambda$ and the stationary points $p_\lambda$ and $q_\lambda$. The stable and unstable manifolds of all the chain-recurrent sets intersect transversally, so the flow is Morse-Smale on $M$.

That the bifurcation described takes place for $X$ in an open subset of $\Psi^{sn}$ follows since **X1** - **X7** each hold in open sets in $\Psi^{sn}$ and the intersection of these subsets is nonempty. All versal deformations are included in condition **X3** and so the result holds for arbitrary versal deformations. $\square$

**Remark 5.3** It is clear that the above procedure may be applied to systems with more complicated pleatings than that described. In fact, one may generalise Definition 3.3 to include systems which 'enclose' branches of the stable or unstable sets, as in the G-horseshoe, without affecting the results. The only sort of behaviour which will ruin the theory is when a fixed point splits $S^s(p)$ (or $S^u(p)$) into different strips which may have non-incremental twist signatures, thereby invalidating Theorem 4.3.

**ACKNOWLEDGMENTS**

This work has been sponsored by the National Science Foundation (RG), the University of Texas at Austin (RG), Northwestern University (TY) and the University of Chicago (TY). The authors especially wish to thank Alec Norton for encouragement.

**AMS Subject Classification:** 34C23, 58F09, 57M25.